\documentclass[]{interact}

\usepackage{epstopdf}
\usepackage[caption=false]{subfig}

\usepackage[numbers,sort&compress]{natbib}
\usepackage{color}
\bibpunct[, ]{[}{]}{,}{n}{,}{,}
\renewcommand\bibfont{\fontsize{10}{12}\selectfont}
\makeatletter
\def\NAT@def@citea{\def\@citea{\NAT@separator}}
\makeatother

\theoremstyle{plain}
\newtheorem{theorem}{Theorem}[section]
\newtheorem{lemma}[theorem]{Lemma}

\newtheorem{proposition}[theorem]{Proposition}

\theoremstyle{definition}
\newtheorem{definition}[theorem]{Definition}
\newtheorem{example}[theorem]{Example}

\theoremstyle{remark}

\newcommand{\R}{{\mathbb R}}
\newcommand{\disp}{\displaystyle}
\newcommand{\abs}[1]{\left\vert#1\mathstrut\right\vert}

\begin{document}


\title{Necessary  conditions for optimal control problems with sweeping systems and end point constraints}

\author{
\name{M.~d.~R. de Pinho\textsuperscript{a},  M.~Margarida A.~Ferreira \textsuperscript{a} and Georgi Smirnov\textsuperscript{b}\thanks{CONTACT G.~Smirnov. Email: smirnov@math.uminho.pt}}
\affil{\textsuperscript{a}Universidade do Porto, Faculdade de Engenharaia, DEEC, SYSTEC, Porto, Portugal. \textsuperscript{b}Universidade do Minho, Dep. Matem\'{a}tica, Centro de 
F\'{i}sica, Campus de Gualtar, Braga, Portugal.}
}

\maketitle

\begin{abstract}
We generalize the Maximum Principle for free end point optimal control problems involving sweeping systems derived in \cite{nosso_2019} to cover the case where the end point is constrained to take values in a certain set. As in \cite{nosso_2019}, an ingenious  smooth approximating family of standard differential equations plays a crucial role.
\end{abstract}

\begin{keywords}
Sweeping systems, optimal control, necessary conditions.
\end{keywords}

\section{Introduction}
Sweeping  processes are evolution differential inclusions involving the normal cone to a set.  They were introduced  in the seminal paper  \cite{Mo74}  by J.J. Moreau in  the context of  plasticity and friction theory. Since then, there has been an increasing interest in sweeping systems with its range of application covering now problems from mechanics, engineering, economics  and  crowd motion problems; see, for example,   \cite{Addy}, \cite{MoCa17}, \cite{KuMa00}, \cite{Maury},  \cite{Thib2016} and \cite{nosso_2019}.

In recent years, there has been considerable research on optimal control problem involving controlled sweeping systems of the form
\begin{equation}\label{SP}
\dot x(t) \in f(t,x(t),u(t))- N_{C(t)}(x(t)), ~u(t)\in U, ~~x(0) \in C_0.
\end{equation}
In this respect, we refer the reader to, for example,
 \cite{ArCo17},  \cite{BrKr},  \cite{MoCa17}, \cite{CoPa16},    \cite{CoHeHoMo}, \cite{KuMa00}, \cite{zeidan2020}, and \cite{nosso_2019} (see also accompanying correction \cite{correction_2019}). 
 A remarkable aspect of \eqref{SP} is that the presence of the  normal cone in this dynamics destroys the  regularity  of the corresponding  differential inclusion under which  classical results on  differential inclusions and optimal control have been built.

Assuming that the set $C$ in \eqref{SP} is time independent, necessary conditions in the form of a maximum principle  for  optimal control  problems involving such systems are derived in  \cite{ArCo_IFAC17}, \cite{BrKr} and \cite{nosso_2019}.
A special feature of  \cite{nosso_2019} is that it relies on an approximating sequence  of optimal control problems differing from the original problem insofar as  \eqref{SP} is replaced  by a   differential equation of the form  
\begin{equation}\label{Sgamma}\dot x_{\gamma_k}(t)=  f(t,x_{\gamma_k}(t),u(t))-\gamma_k e^{\gamma_k \psi(x_{\gamma_k}(t))}\nabla \psi(x_{\gamma_k}(t)) \end{equation}
for some positive sequence $\gamma_k\to +\infty$. Later, similar techniques have been applied  to more general problems in \cite{zeidan2020}.

In this paper we generalize  the Maximum Principle proved in \cite{nosso_2019}   to cover problems with additional end point constraints. Our problem of interest is
$$
(P)
\left\{
\begin{array}{l}
\mbox{Minimize } \; 
\phi(x(T))\\[2mm]
\mbox{over processes $(x,u)$ such that }\\[2mm]
\hspace{8mm} \dot x(t) \in f(t,x(t),u(t))- N_{C(t)}(x(t)), \hspace{0.2cm}\mbox{a.e.}\ \ t\in[0,T],\\[2mm]
\hspace{8mm} u(t)\in U, \ \ \;\, \mbox{a.e.}\ \ t\in[0,T],\\[2mm]
\hspace{8mm} (x(0),x(T)) \in C_0\times C_T { \subset C(0)\times C(T)},
\hspace{2mm}\end{array}
\right.
$$
where  $T>0$ is fixed, $\phi:[0,T]\times\R^n\to \R$, $f:[0,T]\times\R^n\times \R^m\to \R^n$, $U \subset \R^m$ 
 and 
\begin{equation}
\label{set:C}
C(t):=\left\lbrace x\in \R^n: ~\psi(t,x)\leq 0\right\rbrace
\end{equation}
for some  function $\psi:[0,T]\times \R^n\to \R$.
A preliminary attempt to generalize the results of \cite{nosso_2019} to cover problems with end point constraints is submitted to the IEEE CDC Conference 2021 \cite{submetido}. This paper however greatly diverges from  \cite{submetido}; for example, here we assume the set $C$ to be  \emph{time dependent} whereas in \cite{submetido} it is a constant set.

\section{Preliminaries}

In this section, we introduce a summary of the notation and  state  the assumptions  on the data of $(P)$ enforced throughout.  Furthermore, we extract information from the assumptions defining extra functions and establishing relations  crucial for  the forthcoming analysis.
\medskip

\noindent {\bf Notation}

For a set $S\subset \R^n$, $\partial S$, $\text{cl}S$ and $\text{int}\, S$ denote the \emph{boundary}, \emph{closure}   and   \emph{interior} of $S$. 

If  $g:\R^p\to \R^q$, $\nabla g$ represents the derivative and $\nabla^2g$ the second derivative.   
If $g:\R\times \R^p\to \R^q$,  then $\nabla_x g$ represents the derivative w.r.t. $x\in \R^p$ and $\nabla^2_xg$ the second derivative, while $\partial_t g(t,x)$ represents the derivative w.r.t. $t\in \R$.   

The Euclidean norm or the induced matrix
norm on $\R^{p\times q}$ is denoted by $\abs{\,\cdot\,}$.  We denote by $B_n$ the closed unit ball in $\R^n$  centered at the origin. For some $A\subset \R^n$, $d(x,A)$ denotes the distance between $x$ and $A$. 

The space  $L^{\infty}([a,b];\R^p)$ (or simply $L^{\infty}$ when the domains are clearly understood) is the Lebesgue space of essentially 
bounded functions  $h:[a,b]\to\R^p$. We say that $h\in BV([a,b];\R^p)$ if $h$ is a function of bounded variation. The space of continuous functions is denoted by $C([a,b];\R^p)$.

Standard concepts from nonsmooth analysis will also be used. Those can be found in  \cite{Cla83}, \cite{Mor05} or \cite{Vinter}, to name but a few. 
The \emph{Mordukhovich} normal cone to a set $S$  at $s\in S$  is denoted by 
$N_{S}^L(s)$   and 
 $\partial^L f(s)$ is the
\emph{Mordukhovich} subdifferential of $f$ at $s$ (also known as \emph{limiting subdifferential}). 

For any set $A\subset \R^n$, $\text{cone}\: A$ is the cone generated by the set $A$.
\bigskip

We now turn to problem $(P)$.  We first introduce  the definition of admissible processes.

\begin{definition}
A pair  $(x,u)$ is called an admissible process for $(P)$ when  $x$ is an absolutely continuous function and  $u$ is a measurable function  satisfying the constraints of $(P)$.
\end{definition}

\noindent {\bf Assumptions on the data of $\mathbf{(P)}$}
\begin{itemize}
\item[A1:] The function $f$ is  continuous, $x\to f(t,x,u)$  is continuously differentiable  for all $(t,u)\in [0,T]\times \R^m$  and there exists a $M>0$ such that $|f(t,x,u)|\leq M$ and $|\nabla_x f(t,x,u)|\leq M$ for all $(t,x,u)\in [0,T]\times \R^n\times U$.
\item[A2:]  For each $(t,x)$,  the set $f(t,x,U)$ is  convex. 

\item[A3:] The function $\psi$ is  $C^2$  and,  for any $t\in [0,T]$, we have
\begin{equation}
 \disp \lim_{|x| \rightarrow \infty} \psi(t,x) = + \infty,\label{psi:pdef}
 \end{equation} 
 and
there exist constants $\beta>0$ and $\eta >0$  such that
\begin{equation}
\psi(t,x) \geq -\beta \Longrightarrow |\nabla_x \psi (t,x) | > \eta  \text{ for all } (t,x)\in [0,T]\times  \R^n. \label{psi:bder}
\end{equation}
\item[A4:] The set $U$ is compact.
\item[A5:] The sets  $C_0$ and $C_T$ are   closed.
\item[A6:] There exists a constant $L_\phi$ such that $|\phi(x)-\phi(x')|\leq L_{\phi}|x-x'|$ for all $x, x' \in \R^n$.
\end{itemize}
 
The continuity of $\psi$ together with \eqref{psi:pdef} ensure that $C(t)$ is closed and bounded  for all $t\in [0,T]$. Consequently, $C_0$ and $C_T$ are also both compact sets. 
\bigskip

Let $x(\cdot)$ be a solution to the differential inclusion:
$$ \dot x(t) \in f(t,x(t),U)- N_{C(t)}(x(t)).$$
Under our assumptions, measurable selection theorems assert the existence of measurable functions $u$ and $\xi$ such that 
$u(t) \in U$, $\xi(t) \in N_{C(t)}(x(t))$ a.e.  $t\in [0,T]$ and 
$$\dot x(t)= f(t,x(t),u(t))-\xi(t)\nabla_x\psi(t,x(t))\text{ a.e.} t\in [0,T].$$
Let $\mu$ be such that
$$
\max\left\{(|\nabla_x\psi(t,x)||f(t,x,u)|+|\partial_t\psi(t,x)|)+1: ~t\in [0,T], u\in U, x\in C(t)+B_n\right\}\leq\mu.
$$
Observe that the supremum is indeed finite.

For any $t$ such that $\psi (t,x(t))=0$ and $\dot{x}(t)$ exists, we have
$$\begin{array}{c}
0=\frac{d}{dt}\psi(t,x(t))=\langle \nabla_x\psi(t,x(t)),\dot{x}(t)\rangle +\partial_t\psi(t,x(t))
\\[2mm]
=\langle \nabla_x\psi(t,x(t)),f(t,x(t),u(t))\rangle-\xi(t)| \nabla_x\psi(t,x(t))|^2 +\partial_t\psi(t,x(t))
\end{array}$$
Hence
$$
\xi(t)=\frac{1}{| \nabla_x\psi(t,x(t))|^2}(\langle \nabla_x\psi(t,x(t)),f(t,x(t),u(t))\rangle+\partial_t\psi(t,x(t)))
\leq \frac{\mu}{\eta^2}.
$$

Define the function 
$$\mu(\gamma)=\frac{1}{\gamma}\log\frac{\mu}{\eta^2\gamma}.$$
Consider  now a sequence $\{\sigma_k\}$ such that $\sigma_k\downarrow 0$. Let $\{\gamma_k\}$  be a sequence such that
$$C(t)\subset \text{int}\: C^k(t)= \text{int}\: \left\{x: \psi(t,x)-\sigma_k\leq \mu_k\right\}$$
where 
$$\mu_k=\mu(\gamma_k).$$
Consider   $x_k$ to be a solution to the differential equation
$$
\dot{x}_k(t)=f(t,x_k(t),u_k(t))-\gamma_k e^{\gamma_k(\psi(t,x_k(t))-\sigma_k)}\nabla_x\psi(t,x_k(t))
$$
for some $u_k(t) \in U$ a.e. $t\in [0,T]$.
Take any  $t\in [0,T]$ such that $\dot{x}_k(t)$ exists and $\psi (t,x_k(t))-\sigma_k=\mu_k$. We then have
$$ \begin{array}{c}
\frac{d}{dt}\psi(t,x_k(t))\\[2mm]=\langle \nabla_x\psi(t,x_k(t)),f(t,x_k(t),u_k(t))\rangle-\gamma_ke^{\gamma_k(\psi(t,x_k(t))-\sigma_k)}|\nabla_x\psi(t,x_k(t))|^2+
\partial_t\psi(t,x_k(t))
\\[2mm]
\leq \mu-1 -\eta^2\gamma_k e^{\gamma_k\mu_k}=-1.
\end{array}$$
Thus,  if $x_k(0)\in C^k(0)$, then $x_k(t)\in C^k(t)$ for all $t\in [0,T]$, and
\begin{equation}
\gamma_ke^{\gamma_k(\psi(t,x_k(t))-\sigma_k)}\leq \gamma_ke^{\gamma_k\mu_k}= \frac{\mu}{\eta^2}.
\label{star}
\end{equation}
It follows that, for all $k$, we have
$$
|\dot{x}_k(t)|\leq ({\rm const}).
$$

We are now in position to state our first Theorem. This is akin to Theorem 4.1 in \cite{zeidan2020} (see also Lemma 1 in \cite{nosso_2019} when $\psi$ is independent of $t$ and convex) but  we now consider a different approximating sequence of control systems. This guarantees the  estimation \eqref{star},   greatly simplifying  the proof of Theorem \ref{T1}. 
\begin{theorem}
\label{T1}
Let $\{x_k\}$ be a sequence of solutions of Cauchy problems
\begin{equation}
\label{e1}
\dot{x}_k(t)=f(t,x_k(t),u_k(t))-\gamma_k e^{\gamma_k(\psi(t,x_k(t))-\sigma_k)}\nabla_x\psi(t,x_k(t)),\;\;\; x_k(0)=b_k\in C^k(0).
\end{equation}
If $b_k\rightarrow x_0$, then there exists a subsequence $\{x_k\}$ (we do not relabel) converging uniformly  to $x$, a unique solution
to the Cauchy problem
\begin{equation}
\label{e2}
\dot{x}(t)\in f(t,x(t),u(t))-N_{C(t)}(x(t)),\;\;\; x(0)=x_0,
\end{equation}
where $u$ is a  measurable function such that $u(t)\in U$ a.e. $t\in [0,T]$.

If, moreover, all the controls $u_k$ are equal, i.e., $u_k=u$, then the subsequence converges to a unique solution of (\ref{e2}), i.e., any solution of 
\begin{equation}
\label{e3}
\dot{x}(t)\in f(t,x(t),U)-N_{C(t)}(x(t)),\;\;\; x(0)=x_0\in C(0)
\end{equation}
can be approximated by solutions of (\ref{e1}).
\end{theorem}

\begin{proof}
Consider the sequence $\{x_k\}$, where $x_k$ solves \eqref{e1}. Recall that  $x_k(t)\in C^k(t)$ for all $t\in[0,T]$, and 
$$
|\dot{x}_k(t)|\leq ({\rm const})\;\;\;{\rm and}\;\;\; \xi_k(t)=\gamma_k e^{\gamma_k(\psi(t,x_k(t))-\sigma_k)}\leq ({\rm const}).
$$
Then there exist subsequences (we do not relabel) weakly-$*$ converging in $L^{\infty}$ to some  $v$ and $\xi$.
Hence
$$x_{k}(t)=x_0+\displaystyle\int_0^t \dot x_{k}(s) ds \longrightarrow x(t)=x_0+
\displaystyle\int_0^t  v(s)ds, ~\forall ~t\in [0,1],
$$
for an absolutely continuous function $x$. Obviously, $x(t)\in C(t)$ for all $t\in [0,T]$.

We have
\begin{equation}
\dot{x}_k(t)\in f(t,x_k(t),U)-\xi_k(t)\nabla_x \psi(t,x_k(t)).\label{ap1}
\end{equation}

Let $A\subset R^n$ and $ z\in R^n$. We denote the support function of $A$ at $z$ by $S(z,A)=\sup\{\langle z,a\rangle\mid a\in A\}$.
Inclusion (\ref{ap1}) is equivalent to the condition
$$
\langle z, \dot{x}_k(t)\rangle\leq S(z,f(t,x_k(t),U))-\xi_k(t)\langle z,
 \nabla_x\psi(t,x_k(t)),\;\;\;\forall\: z\in R^n.
 $$
Integrating this inequality, we get
 $$
 \left\langle z,\frac{x_k(t+\tau)-x_k(t)}{\tau}\right\rangle\leq \frac{1}{\tau}\int_t^{t+\tau}(S(z,f(s,x_k(s),U))-\xi_k(s)\langle z,
 \nabla_x\psi(s,x_k(s))\rangle )ds
 $$
 $$
 =\frac{1}{\tau}\int_t^{t+\tau}(S(z,f(s,x_k(s),U))-\xi_k(s)\langle z, \nabla_x\psi(s,x(s))\rangle 
 $$
 $$
 +
\xi_k(s)\langle z, \nabla_x\psi(s,x(s))- \nabla_x\psi(s,x_k(s))\rangle )ds.
$$
 Passing to the limit as $k\rightarrow\infty$, we obtain
  $$
 \left\langle z,\frac{x(t+\tau)-x(t)}{\tau}\right\rangle\leq \frac{1}{\tau}\int_t^{t+\tau}(S(z,f(s,x(s),U))-\xi(s)\langle z,
 \nabla_x\psi(s,x(s))\rangle )ds.
  $$
Let $t\in[0,1]$ be a Lebesgue point of  $x$ and $\xi$. Passing to the limit as $\tau\downarrow 0$, we have
$$
\langle z,\dot{x}(t)\rangle\leq S(z,f(t,x(t),U))-\xi(t)\langle z,
 \nabla_x\psi(t,x(t))\rangle.
 $$
Since $z\in R^n$ is an arbitrary vector and the set $f(t,x(t),U)$ is convex, we conclude that
$$
\dot{x}(t)\in f(t,x(t),U)-\xi(t)\nabla_x\psi(t,x(t)).
$$
By the Filippov lemma there exists a measurable control $u(t)\in U$ such that
$$
\dot{x}(t)= f(t,x(t),u(t))-\xi(t)\nabla_x\psi(t,x(t)).
$$
Observe  that $\xi$ is zero if $\psi(t,x(t))<0$.

If $u_k=u$ for all $k$, then the sequence $x_k$ converges to the solution of
$$
\dot{x}(t)= f(t,x(t),u(t))-\xi(t)\nabla_x\psi(t,x(t)).
$$
Indeed, to see this it suffices to pass to the limit as $k\rightarrow\infty$ and then as $\tau\downarrow 0$ in the equality
$$
\frac{x_k(t+\tau)-x_k(t)}{\tau}= \frac{1}{\tau}\int_t^{t+\tau}(f(s,x_k(s),u(t))-\xi_k(s) \nabla_x\psi(s,x_k(s)) )ds.
$$

\medskip

We now prove the uniqueness of the solution. We follow the proof of Theorem 4.1 in \cite{zeidan2020}. Notice however that we now consider a special case and not the general case treated  in \cite{zeidan2020}. 

 Suppose that there exist two different solutions of (\ref{e2}): $x_1$ and $x_2$. We have
$$
\frac{1}{2}\frac{d}{dt}|x_1(t)-x_2(t)|^2=\langle x_1(t)-x_2(t),\dot{x}_1(t)-\dot{x}_2(t)\rangle
$$
$$
=\langle x_1(t)-x_2(t),f(t,x_1(t),u(t))-f(t,x_2(t),u(t))\rangle
$$
$$
-\langle x_1(t)-x_2(t),\xi_1(t)\nabla\psi(t,x_1(t))-\xi_2(t)\nabla\psi(t,x_2(t))\rangle.
$$
If $\psi(t,x_1(t))<0$ and $\psi(t,x_2(t))<0$, then $\xi_1(t)=\xi_2(t)=0$ and we obtain
$$
\frac{1}{2}\frac{d}{dt}|x_1(t)-x_2(t)|^2\leq L_f|x_1(t)-x_2(t)|^2.
$$
Suppose that $\psi(t,x_1(t))=0$, then by the Taylor formula we get
$$
\psi(t,x_2(t))=\psi(t,x_1(t))+\langle \nabla_x\psi(t,x_1(t)),x_2(t)-x_1(t)\rangle
$$
$$
+\frac{1}{2}\langle x_2(t)-x_1(t), \nabla_x^2\psi(t,\theta x_2(t)+(1-\theta )x_1(t))( x_2(t)-x_1(t))\rangle,
$$
where $\theta\in [0,1]$. Since $\psi(t,x_2(t))\leq  0$, we have
$$
\langle \nabla_x\psi(t,x_1(t)),x_2(t)-x_1(t)\rangle
$$
$$
\leq -
\frac{1}{2}\langle x_2(t)-x_1(t), \nabla_x^2\psi(t,\theta x_2(t)+(1-\theta )x_1(t))( x_2(t)-x_1(t))\rangle\leq ({\rm const})
|x_1(t)-x_2(t)|^2.
$$
In the same way 
$$
\langle \nabla_x\psi(t,x_2(t)),x_1(t)-x_2(t)\rangle\leq ({\rm const})|x_1(t)-x_2(t)|^2.
$$
Thus we have
$$
\frac{1}{2}\frac{d}{dt}|x_1(t)-x_2(t)|^2\leq ({\rm const})|x_1(t)-x_2(t)|^2.
$$
Hence $|x_1(t)-x_2(t)|=0$. 
\end{proof}

\section{Approximating Family of Optimal Control Problems}
In this section we define an approximating family of optimal control problems to $(P)$ and we state the corresponding necessary conditions. 

Let $(\hat{x},\hat{u})$  
be a global solution to $(P)$ and consider a sequence $\{\gamma_k\}$ as defined above.
Let  $\hat  x_k(\cdot)$ be the solution to
\begin{equation}\label{starstar}\left\{\begin{array}{rcl}
\dot{x}(t)& = & f(t,x(t),\hat u(t))-\gamma_k e^{\gamma_k(\psi(t,x(t))-\sigma_k)}\nabla_x\psi(t,x(t)),\\[2mm]
x(0)&= & \hat x(0).
\end{array}
\right.
\end{equation}
Set   $\epsilon_k=|\hat x_k(T)-\hat x(T)|$. It follows from  Theorem \ref{T1}  that $\epsilon_k\downarrow 0$. 
Take $\alpha>0$ and define the problem
$$
(Q_k^\alpha)
\left\{
\begin{array}{l}
\mbox{Minimize } \; 
\phi(x(T))+|x(0)-\hat{x}(0)|^2+\alpha\disp\int_0^T|u(t)-\hat{u}(t)|dt\\[2mm]
\mbox{over processes $(x,u)$ such that }\\[2mm]
\hspace{8mm} \dot x(t) = f(t,x(t),u(t))-\nabla_x e^{\gamma_k(\psi(t,x(t))-\sigma_k)} \hspace{0.2cm}\mbox{a.e.}\ \ t\in[0,T],\\[2mm]
\hspace{8mm} u(t)\in U~~~~\mbox{a.e.}\ \ t\in[0,T],\\[2mm]
\hspace{8mm} x(0)\in C_0,~~ x(T)\in C_T+\epsilon_k B_n,
\hspace{2mm}\end{array}
\right.
$$
Clearly, the problem 
$(Q_k^\alpha)$ has admissible solutions.

Consider  the metric space $$W=\{ (c,u)\mid c\in C_0,\; u\in L^{\infty} \text{ with } u(t)\in U\}$$
and  the distance
$$
d_{W}((c_1,u_1),(c_2,u_2))=|c_1-c_2|+\int_0^T|u_1(t)-u_2(t)|dt.
$$

Endowed with $d_{W_k}$,  $W_k$ is a complete metric space.  Take any $(c,u)\in W$ and a solution $y$ to the Cauchy problem
$$\left\{\begin{array}{rcl}
 \dot y(t) & = &  f(t,y(t),u(t))-\nabla_x e^{\gamma_k(\psi(t,y(t))-\sigma_k)} \hspace{0.2cm}\mbox{a.e.}\ \ t\in[0,T],\\[2mm]
 y(0) & = & c.
 \end{array}
 \right.
 $$
Under our conditions,
the function 
$$(c,u)~\rightarrow ~  \phi(y(T))+ | c - \hat{x}(0) |^2+\alpha \int_{0}^{T} 
| u-\hat{u}|~dt$$  is continuous on $(W,d_{W})$ and bounded below.

 Appealing to Ekeland's Theorem we deduce the existence of a pair $(x_k,u_k)$ solving the following problem
$$
(AQ_k)
\left\{
\begin{array}{l}
\mbox{Minimize } \; 
\Phi(x,{u})= \phi(x(T))+|x(0)-\hat{x}(0)|^2+\alpha\disp\int_0^T|u(t)-\hat{u}(t)|dt\\[2mm]
\qquad\qquad +\epsilon_k\left(|x(0)-x_k(0)|+\disp\int_0^T|u(t)-u_k(t)|dt\right),\\[2mm]
\mbox{over processes $(x,u)$ such that }\\[2mm]
\hspace{8mm} \dot x(t) = f(t,x(t),u(t))-\nabla_x e^{\gamma_k(\psi(t,x(t))-\sigma_k)} \hspace{0.2cm}\mbox{a.e.}\ \ t\in[0,T],\\[2mm]
\hspace{8mm} u(t)\in U~~~~\mbox{a.e.}\ \ t\in[0,T],\\[2mm]
\hspace{8mm} x(0)\in C_0,~~ x(T)\in C_T+\epsilon_k B_n,
\hspace{2mm}\end{array}
\right.
$$
\begin{lemma}\label{lemma:approx}
Take $\gamma_k\to \infty$ and  $\epsilon_k \to 0$ as defined above. For each $k$, let  $(x_k,u_k)$ be the solution to $(AQ_k)$. Then there exists a  subsequence (we do not relabel) such that 
$$u_k(t)\rightarrow \hat u(t)~\text{a.e.}, \quad x_k\rightarrow \hat x \text{ uniformly in } [0,1].$$
\end{lemma}
\begin{proof}

Take any admissible process  $(x,u)$ for $(P)$. Then there exists a $K>0$ such that 
$$
 |x(0)-x_k(0)|+\int_0^T|u(t)-u_k(t)|dt\leq K
 $$
Let us assume that there exists a  $\delta>0$  such that
\begin{equation}\label{contradiction}
|x_k(0)-\hat{x}(0)|^2+\alpha\int_0^T |u_k(t)-\hat{u}(t)|dt\geq\delta,
\end{equation}
for all $k$. 
It is an easy task to see that, for large $k$,  we have $\delta>8\epsilon_kK$. 

\medskip

Now, we deduce  from Theorem \ref{T1} that $\{x_k\}$ uniformly converges to an admissible  solution $\tilde{x}$ to  
$(P)$.  
Taking into account the optimality of $(\hat x,\hat u)$ and \eqref{contradiction},  we conclude that
\begin{equation}
\label{ee1}
\Phi (x_k,u_k)\geq \phi(x_k(T))+\delta\geq \phi(\tilde x(T))+\delta/2\geq \phi(\hat{x}(T))+\delta/2.
\end{equation}
Again from Theorem \ref{T1},  it is also an easy task to see that $\hat x_k$  (see \eqref{starstar}) satisfies the inequality
\begin{equation}
\label{ee2}
\phi(\hat x_k(T))\leq\phi(\hat{x}(T))+\delta/8,
\end{equation}
 for   $k$ is sufficiently large; in this respect, notice that each $\hat x_k$ corresponds  to the control $\hat u$.

It follows from  the optimality of $(x_k,u_k)$ and \eqref{ee1} that
\begin{equation}
\label{e}
\Phi(\hat x_k,\hat{u})\geq\Phi(x_k,u_k)\geq \phi(\hat{x}(T))+\delta/2.
\end{equation}
On the other hand, from (\ref{ee2}), we get 
$$\begin{array}{c}
\Phi(\hat x_k,\hat{u})=\\[2mm]
\phi(\hat x_k(T))+|\hat x_k  (0)-\hat x(0)|^2 +\epsilon_k\left( |\hat x_k(0)-x_k(0)|+\disp\int_0^T|\hat u(t)-u_k(t)|dt\right)\\[2mm]
\leq\phi(\hat{x}(T))+\delta/8+|\hat x_k(0)-\hat{x}(0)|^2+\epsilon_kK\\[2mm]<\phi(\hat{x}(T))+\delta/4,
\end{array}
$$
for sufficiently large $k$, contradicting \eqref{e}. Then we conclude that
$$
|x_k(0)-\hat{x}(0)|^2+\alpha \int_0^1 |u_k(t)-\hat{u}(t)|dt\rightarrow 0,\;\; k\rightarrow \infty
$$ 
and the result follows.
\end{proof}

We now finish this section with  the statement of the optimality necessary conditions for the family of problems $(AQ_k)$. These can be seen as a  direct consequence of  Theorem 6.2.1 in \cite{Vinter}, for example.
\begin{proposition}\label{lemma:nco}
For each k, let $(x_k,u_k)$ be a solution to $(AQ_k)$. Then there exist absolutely continous functions $p_k$ and scalars  $\lambda_k\geq 0$ such that
\begin{itemize}
\item[{\bf(a)}] (nontriviality condition) \begin{equation}
\label{not0}
  \lambda_k+|p_{k}(T)|  =1,
\end{equation}

\item[{\bf(b)}] (adjoint equation)
\begin{equation}
\label{dynamic}
\dot{p}_{k} = -(\nabla_x f_{k})^* p_{k} 
+\gamma_k e^{\gamma_k (\psi_{k}-\sigma_k)}\nabla^2_x\psi_{k}p_{k}
+\gamma_k^2e^{\gamma_k (\psi_{k}-\sigma_k)}\nabla_x\psi_{k}\langle\nabla_x\psi_{k},p_{k}\rangle, 
\end{equation}
where the superscript $*$ stands for transpose,

\item [{\bf(c)}] (maximization condition)
\begin{equation}
\label{max}
\disp \max_{u\in U}\left\{ \langle f(t,x_{k}, u)  ,  p_{k} \rangle -  \alpha \lambda_k|u-\hat{u}|
-\epsilon_k \lambda_k|u-u_k|\right\}
\end{equation}
is attained at  $u_k (t) $, for almost every $t\in [0,T]$,
\item [{\bf(d)}] (transversality condition)
\begin{eqnarray}
 ( p_{k}(0), - p_{k}(T)) \hspace{5cm}\nonumber \\[2mm] 
 \in \lambda_k\left(2(x_k(0)-\hat{x}(0)),  \partial^L  \phi (x_{k}(T))\right)
 +  N_{C_0}^L(x_{k}(0))\times N_{C_T+\epsilon_kB_n}^L(x_{k}(T)). \label{trans}
\end{eqnarray}
\end{itemize}
\end{proposition}

To simplify the notation above, we drop the $t$ dependance in $p_k$, $\dot p_k$, $x_k$, $u_k$, $\hat x$ and $\hat u$. Moreover,  in (b), we write $\psi_k$ instead of $\psi(t,x_k(t))$, $f_k$ instead of $f(t,x_k(t),u_k(t))$. The same holds for the derivatives of $\psi$ and $f$.

\section{Maximum Principle for $(P)$}

In this section we establish our main result, a Maximum Principle for $(P)$.  This is done by taking  limits of the conclusions of Proposition \ref{lemma:nco}, following closely the analysis done in the proof of \cite[Theorem 2]{nosso_2019}. However, some changes are called for since here $C$ is time dependent, while  in \cite{nosso_2019} $C$ is constant. 

First, observe that
$$\begin{array}{c}
\frac{1}{2} \frac{d}{dt} |p_k(t)|^2 = \\[2mm]
- \langle \nabla_x f_k p_k , p_k \rangle +  
\gamma_k e^{\gamma_k(\psi_k-\sigma_k)} \langle \nabla_x^2\psi_kp_k, p_k \rangle +  \gamma_k^2e^{\gamma_k(\psi_k-\sigma_k)}\left(\langle \nabla_x\psi_k, p_k \rangle \right)^2 \\[2mm]
\geq - \langle \nabla_x f_k p_k , p_k \rangle +  
\gamma_k e^{\gamma_k(\psi_k-\sigma_k)} \langle \nabla_x^2\psi_kp_k, p_k \rangle\\[2mm]
 \geq \ - M | p_k|^2 +\gamma_k e^{\gamma_k(\psi_k-\sigma_k)} \langle \nabla_x^2\psi_kp_k, p_k \rangle,
\end{array}$$
where $M$ is the constant of (A1). Taking into account hypothesis  (A3) and \eqref{star} we deduce the existence of a constant $K_0>0$ such that
$$\frac{1}{2} \frac{d}{dt} |p_k(t)|^2\geq  -K_0| p_k|^2.$$
This last inequality leads to 
$$
| p_k(t)|^2 \ \leq \ e^{2 K_0 (T-t)} | p_k(T)|^2 \leq  \ e^{2 K_0T} |p_k(T)|^2.
$$
Since, by (a) of Proposition \ref{lemma:nco},  $|p_k(T)|\leq 1$, we deduce from the above that
there exists  $M_0>0$ such that
\begin{equation}\label{pk1}
| p_k(t)| \ \leq  M_0.
\end{equation}

Next, we prove that the sequence $\{\dot p_k\}$ is uniformly bounded in $L^1$. 
We start by establishing some inequalities that differ from analogous in \cite{nosso_2019} due, among other things,  to the $\psi$ dependence on $t$.

We have
$$\begin{array}{c}
\frac{d}{dt} \left| \langle  \nabla_x\psi_ k, p_ k \rangle \right| = 
\left( \langle \nabla^2_x\psi_ k \dot{x}_ k, p_ k \rangle
+\langle\partial_t\nabla_x\psi_k,p_k\rangle
+ \langle\nabla_x\psi_ k,\dot{p}_ k\rangle \right)\,  {\rm sign}\left( \langle  \nabla_x\psi_ k, p_ k \rangle \right)\\[2mm]
=
\left( \langle  p_ k, \nabla^2_x\psi_ k f_ k \rangle
-  \gamma_k e^{\gamma_k\psi_ k(\psi_ k-\sigma_k)}  \langle p_ k, \nabla^2\psi_ k   \nabla_x\psi_ k  \rangle \right. \\[2mm]
+\langle\partial_t\nabla_x\psi_k,p_k\rangle
- \langle  \nabla_x \psi_ k, (\nabla_x f_ k)^* p_ k \rangle 
+ 
\gamma_k e^{\gamma_k(\psi_ k-\sigma_k)}  \langle  \nabla_x\psi_ k ,  \nabla^2_x\psi_ k  p_ k\rangle
 \\[2mm]
+\left.  \gamma_k^2e^{\gamma_k(\psi_ k-\sigma_k)} | \nabla_x\psi_ k|^2 \langle\nabla_x\psi_ k,p_ k\rangle \right) \, {\rm sign}\left( \langle  \nabla_x\psi_ k, p_ k \rangle \right).
\end{array}
$$
Moreover, we also have
\begin{eqnarray*}
& \gamma_k^2  \disp\int_{0}^{T} e^{\gamma_k(\psi_ k-\sigma_k)} | \nabla_x\psi_ k|^2  \left| \langle \nabla_x\psi_ k,p_ k\rangle \right| \ dt  \\[2mm]
&
=\left| \langle \nabla_x\psi(T,x_ k(T)),p_ k(T) \rangle \right|  - \left| \langle \nabla_x\psi(0,x_ k(0)),p_ k(0) \rangle \right|
\\[2mm]
& \quad \quad +  \disp \int_{0}^{T} \left( \langle  \nabla_x \psi_ k, (\nabla_x f_ k)^* p_ k \rangle 
-\langle\partial_t\nabla_x\psi_k,p_k\rangle - \langle  p_ k, \nabla_x^2\psi_ k f_ k \rangle \right)  \, {\rm sign}\left( \langle  \nabla_x\psi_ k, p_ k \rangle \right) \ dt\\[2mm]
& \leq \ M_1,
\end{eqnarray*}
for some $M_1>0$.

We now concentrate on  $ \gamma_k^2   \disp\int_{0}^{T} e^{\gamma_k(\psi_ k-\sigma_k)} |\nabla_x\psi_ k|  \left| \langle \nabla_x\psi_ k,p_ k\rangle \right| \ dt $.

Set
$L_k=e^{\gamma_k(\psi_ k-\sigma_k)} |\nabla_x\psi_ k|  \left| \langle \nabla_x\psi_ k,p_ k\rangle\right|$ and observe that
$$  \gamma_k^2   \int_{0}^{T} L_kdt=\gamma_k^2   \int_{\{t: |\nabla_x \psi_ k| < \eta\}} \hspace{-0.7cm}L_kdt+
\gamma_k^2   \int_{\{t: |\nabla_x \psi_ k| \geq \eta\}}\hspace{-0.7cm} L_k dt.$$
Then, from (A3) we deduce that
 $$
 \gamma_k^2   \int_{0}^{T} e^{\gamma_k(\psi_ k-\sigma_k)} |\nabla_x\psi_ k|  \left| \langle \nabla_x\psi_ k,p_ k\rangle \right| \ dt  
 $$
 $$
\leq   \gamma_k^2  e^{-\gamma_k(\beta-\sigma_k)} \eta^2 \max_{t} |p_ k(t)|
+
\gamma_k^2  \int_{\{t: |\nabla_x \psi_ k| \geq \eta\}} \hspace{-1cm}e^{\gamma_k(\psi_ k-\sigma_k)} \frac{|\nabla_x\psi_ k|^2}{ | \nabla_x\psi_ k|} \left| \langle  \nabla_x\psi_ k,p_ k\rangle \right| \ dt 
$$
$$
\leq \gamma_k^2  e^{-\gamma_k (\beta-\sigma_k)} \eta^2 M_0 + \frac{\gamma_k^2 }{\eta}
  \int_{0}^{T} e^{\gamma_k(\psi_ k-\sigma_k)} | \nabla_x\psi_ k|^2 \left| \langle \nabla_x\psi_ k,p_ k\rangle \right|\ dt
  $$
  $$
 \leq \ \eta^2 M_0 + \frac{M_1}{\eta},
  $$
for $k$ large  enough. 

Summarizing, for some $M_2>0$, we have
\begin{equation}
\label{M2}
 \gamma_k^2 \int_{0}^{1} e^{\gamma_k(\psi_ k-\sigma_k)} |\nabla\psi_ k|  \left| \langle \nabla\psi_ k,p_ k\rangle \right| \ dt \  \ \leq \ \eta^2 M_0 + \frac{M_1}{\eta} = M_2.
\end{equation}
Mimicking the analysis conducted in Step 1, b) and c) of  the proof of Theorem 2 in \cite{nosso_2019}
 and taking into account (b) of Proposition \ref{lemma:nco} we conclude that there exist constants $N_1>0$ and $N_2>0$ such that
 \begin{equation}
\label{M8}
\int_0^1 \gamma_k^2e^{\gamma_k(\psi_{\gamma_k}-\sigma_k)} |\langle\nabla\psi_{\gamma_k},p_{\gamma_k}\rangle|dt\leq N_2
\end{equation}
and
\begin{equation}
\label{M7}
\int_{0}^{1} \left| \dot{p}_{\gamma_k}(t)\right|dt \leq N_1, 
\end{equation}
for $k$ sufficiently large. 
With such bounds, we can then mimick the analysis of Step 2  in the proof of Theorem 2 in \cite{nosso_2019} to conclude the existence of functions  $p\in BV([0,1],R^n)$,  $\xi\in L^{\infty}([0,1],R)$, $\xi(t) \geq 0 \ \mbox{ a. e. } t$, $\xi(t) = 0, \ t \in I_b$,  where
$$I_b=\left\{t\in [0,T]:~\psi(t, \hat x(t))<0\right\},$$
 and a finite signed Radon measure $\eta$, null in $ I_b$, such that, for any $z\in C([0,1],R^n)$ 
\begin{equation}\label{Intermideate:a}
\int_0^1 \langle z,dp\rangle=-\int_0^1 \langle z, (\nabla\hat{f})^*p\rangle dt +\int_0^1\xi\langle z,\nabla^2\hat{\psi}p\rangle dt
+\int_0^1\langle z, \nabla\hat{\psi} (t)  \rangle d\eta.
\end{equation}

It is a simple matter to see that there exists a subsequence of  $\{\lambda_k\}$ converging to some $\lambda\geq 0$. This, together with convergence of $p_k$ to $p$ allows us to take limits in (a) and (c) of Proposition \ref{lemma:nco} to deduce that
$$\lambda+|p(T)|=1$$
and
$$\langle p(t), f(t,\hat x(t),u)\rangle-\alpha\lambda |u-\hat{u}(t)| \leq \langle p(t), f(t,\hat x(t),\hat u(t))\rangle\text{   for all } u\in U\text{ and a.e. } t \in [0,T].$$

It remains to take limits of the transversality conditions (d) in Proposition \ref{lemma:nco}.
First, observe that 
$$C_T+\epsilon_kB_n=\left\{x:~d(x,C_T)\leq \epsilon_k\right\}.$$
>From the basic properties of the Mordukhovich normal cone and subdifferential (see \cite{Mor05}, section 1.3.3) we have
$$
N_{C_T+\epsilon_kB_n}^L(x_k(T))\subset \text{ cl cone}\:\partial^L d(x_k(T), C_T)$$
and
$$
N_{C_T}^L(\hat{x}(T))= \text{ cl cone}\:\partial^L d(\hat{x}(T), C_T).
$$
Passing to the limit as $k\to\infty$ we get
$$(p(0),-p(T))\in 
N_{C_0}^L(\hat x(0))\times N_{C_T}^L(\hat x(T))+\{0\}\times \lambda  \partial^L  \phi(\hat x(T)).$$

Now we mimick Step 3 in the proof of Theorem 2 in \cite{nosso_2019}  to remove the dependence of the conditions on the parameter $\alpha$ which is done,  by taking further limits, this time  considering a sequence of $\alpha_i\downarrow 0$. 

We summarize our conclusions in the  following Theorem.
%
%
%
%
%

\begin{theorem}\label{main}
 Let $(\hat{x}, \hat{u})$  be the optimal solution to  $(P)$. 
Suppose that assumption A1--A6 are satisfied. Set $\disp I_b=\left\{t\in [0,T]: ~\psi(t,\hat x(t))<0\right\}$. 

Then there exist a non negative scalar $\lambda$,
$p\in BV([0,T];\R^n)$,   a finite signed Radon measure $\eta$, null in $I_b$, 
$\xi \in L^{\infty}([0,T];\R)$ with
$\disp \xi(t)\geq 0 \text{ a.e. } t, ~\xi(t)=0 \text{ for } t \in I_b$, 
satisfying the following conditions
\begin{itemize}
\item [(a)] 
$\lambda+|p(T)|\neq 0$,

\medskip

\item [(b)] 
for any $z\in C([0,1];\R^n)$ 
$$\begin{array}{c}
\disp\int_0^T \langle z(t),dp(t)\rangle =\\[2mm]
-\disp\int_0^T \langle z(t), ( \nabla_x \hat f(t))^*p(t)\rangle dt+
\disp\int_0^1\!\xi(t) \langle z(t), \nabla^2_x\hat \psi(t) p(t)\rangle dt+\disp\int_0^1\!\langle z(t), \nabla_x \hat \psi(t)\rangle d\eta,\end{array}$$
where the superscript $*$ stands for transpose,  $ \nabla_x \hat f(t) = \nabla_x f(t,\hat x(t),\hat u(t))$, $\nabla_x \hat \psi(t)=\nabla_x \psi(t,\hat x(t))$ and $\nabla^2_x \hat \psi(t)=\nabla^2_x \psi(t,\hat x(t)),$
\medskip

\item [(c)]
$\langle p(t), f(t,\hat x(t),u)\rangle \leq \langle p(t), f(t,\hat x(t),\hat u(t))\rangle$ for  a.e. $t \in [0,T]$
and  all $u\in U,$
\medskip

\item[(d)]
$(p(0),-p(T))\in 
N_{C_0}^L(\hat x(0))\times N_{C_T}^L(\hat x(T))+\{0\}\times \lambda  \partial^L  \phi(\hat x(T))
$
\end{itemize}
\end{theorem}

\section{Example}

We consider two examples.  For the first one, we construct the optimal solution and then  we show that the  necessary conditions given by Theorem \ref{main} are satisfied by all admissible solutions.
The situation changes for the second example where the necessary conditions significantly reduce the set of candidates to the solution.

\begin{example}
Set  $\rho(t)=(1-2t)^2+\frac{1}{4}$ and $C(t)=\{ (x,y)\in R^2\mid x^2+y^2\leq\rho^2(t)\}$. 
Consider now  the control system
$$\left\{\begin{array}{l}
\frac{d}{dt}\left(
\begin{array}{c}
x(t)\\
y(t)
\end{array}
\right)=
\left(
\begin{array}{c}
u(t)\\
0
\end{array}
\right)-N_{C(t)}(x(t),y(t))
,\\[4mm]
u(t)\in [-\mu, 1].
\end{array}
\right.
$$
Observe that $$\left(
\begin{array}{c}
u(t)\\
0
\end{array}
\right)-N_{C(t)}(x(t),y(t))
=
\left(
\begin{array}{c}
u(t)\\
0
\end{array}
\right)
-\lambda(t)\left(
\begin{array}{c}
x(t)\\
y(t)
\end{array}\right)$$
In the above,  $\mu>0$ and $\lambda$ is a $L^1$ function with $\lambda(t)\geq 0$. 

Introducing  polar coordinates
$x=r\cos\phi$, $y=r\sin\phi$, the system takes the form
\begin{eqnarray*}
&& \dot{r}=u\cos\phi-\lambda r,\\
&&\dot{\phi}=-\frac{u\sin\phi}{r}.
\end{eqnarray*}

We consider the following optimal control problem: minimize  
$\Phi(x(T),y(T))=-x(T)$ over the solutions of the above system coupled with $(x(0),y(0))=(0,\frac{\sqrt{3}}{4})$
and $(x(T),y(T))\in C_T$, where  $C_T=\{ (x,y)\in R^2\mid 
x^2+y^2\leq r_T^2\}$. 

Assume that the final time $T$ has the form $T=\frac{1+3\theta}{2}$ for some $\theta$.  We define  $\theta$, $\mu$, and $r_T$  below.

\vspace{5mm}
 The optimal control has the form
$$
\hat{u}(t)=\left\{
\begin{array}{cl}
1, & t\in [0,t_1],\\
1, & t\in [t_1,t_2],\\
-\mu, & t\in ]t_2,T].
\end{array}
\right.
$$
The distance $\hat{r}(t)$ between the origin and the vector $(\hat{x}(t),\hat{y}(t))$ depends on time in the following manner:
$$
\hat{r}(t)=\left\{
\begin{array}{cl}
\sqrt{\hat{y}^2(0)+t^2}, & t\in [0,t_1],\\[2mm]
\rho(t), & t\in [t_1,t_2],\\[2mm]
\sqrt{\hat{y}^2(t_2)+(\hat{x}(t_2)-\mu (t-t_2))^2}, & t\in ]t_2,T].
\end{array}
\right.
$$
So, the trajectory has the following structure: between $t=0$ and $t=t_1=\frac{1}{4}$ the coordinate $x$ grows up,   $\hat{x}(t)=t$, and 
$y$ is constant with 
$\hat{y}(t)\equiv y(0)=\frac{\sqrt{3}}{4}$. At the instant $t_1$ the trajectory reaches the set $C(t_1)$ and then  moves along the
boundary of $C(\cdot)$ for  $t\in [t_1,t_2]$. Finally, at $t=t_2$ it leaves the boundary of $C(\cdot)$ and moves toward $C_T$ along 
the axis $x$ with constant velocity $\dot{\hat{x}}(t)=-\mu$. 
The coordinate $y$ can diminish only if the system moves along the
boundary of $C(\cdot)$. Also, the absolute minimum of $\Phi$ on $C_T$ is achieved at the point $(r_T,0)$. Observe here that  $y(t)$ cannot reach zero
in finite time. Thus, the minimum of $\Phi$, $-x(T)=-\sqrt{r_T^2-y^2(T)}=-\sqrt{r_T^2-y^2(t_2)}$, is achieved if the value of $y(t_2)>0$ is 
as small as possible. Let us now define  $t_2$  as $\frac{1}{2}+\theta$.

\vspace*{5mm}

Now we are ready to define the parameters  $\theta$, $\mu$, and $r_T$. For large $t$ the derivative $\dot{\rho}
(t)$ is so large  that the system cannot move along the boundary of $C(\cdot)$. To determine the moment $t_*$ after which  the system is unable to move on the boundary of $C(\cdot)$, we have to solve the equation
\begin{equation}
\dot{r}(t_*)=1\cdot\cos\phi(t_*)-0\cdot r(t_*)=\dot{\rho}(t_*)=8t_*-4.\label{ex1}
\end{equation}
To find $\phi_*=\phi(t_*)$ we solve the equation
$$
\dot{\phi}=-\frac{u\sin\phi}{r}, \;\;\; t\in [t_1,t_*],
$$
with $\phi(t_1)=\hat{\phi}(t_1)$, $u=1$ and $r(t)=\rho(t)$. It is easy to see that
$$
\phi(t_1)=\phi\left(\frac{1}{4}\right)=\frac{\pi}{3}
$$
and
\begin{equation}
\phi(t)=2\arctan\left(\frac{1}{\sqrt{3}}\exp\left( \arctan ( 2(1-2t))-\frac{\pi}{4}\right)\right),\;\;\; t\in [t_1,t_*].\label{ex2}
\end{equation}

Set $\tau=\tan\frac{\phi_*}{2}$. From (\ref{ex1}) and (\ref{ex2}) we get the system of equations
\begin{eqnarray*}
&& \frac{1-\tau^2}{1+\tau^2}=8t_*-4,\\
&& \tau=\frac{1}{\sqrt{3}}\exp\left( \arctan ( 2(1-2t_*))-\frac{\pi}{4}\right).
\end{eqnarray*} 
This system has a solution $\tau>0$ and 
$$
t_*=\frac{1}{8}\left( 4+ \frac{1-\tau^2}{1+\tau^2}\right) >\frac{1}{2}.
$$
We define $\theta$ as 
$$
\theta=\frac{1}{2}\left(t_*-\frac{1}{2}\right).
$$
In terms of $\theta$ we have
$$
\rho(T)=9\theta^2+\frac{1}{4},\;\;\; \hat{r}(t_2)=\rho(t_2)=4\theta^2+\frac{1}{4}.
$$
Consider the function
$$
r_T^2(\mu)=\hat{y}^2(t_2)+\left(\hat{x}(t_2)-\frac{1}{2}\mu\theta\right)^2=
\hat{r}^2(t_2)-\mu\theta\hat{r}(t_2)\cos\hat{\phi}(t_2)+\frac{\mu^2}{4}\theta^2.
$$
Since $r_T^2(0)=\hat{r}^2(t_2)$ and $\frac{d}{dt}\rho^2(t_2)>0$, for small $\mu>0$ we have
$$
\rho(T)>\hat{r}(t_2)>r_T(\mu)>\hat{r}(t_2)\sin\hat{\phi}(t_2),
$$
and
$$
\Delta=\frac{d}{dt}\rho^2(t_2)-\mu\theta\dot{x}(t_2)+2\mu x(t_2)-\mu^2\theta>0.
$$
We fix such $\mu>0$ and put $r_T=r_T(\mu)$. This choice of $\mu$ guarantees that $r(T)=r_T$, i.e. $(\hat{x}(T),\hat{y}(T))\in C_T$.
Moreover  $\min\Phi=-\hat{x}(T)$. Indeed, consider the trajectory $(x(t),y(t))$, $t\in [t_2,t_2+h]\cup [t_2+h,T]$, of our system, satisfying the following conditions: 
$$\begin{array}{l}(x(t_2),y(t_2))=(\hat{x}(t_2),\hat{y}(t_2)) \text{ and }
x^2(t)+y^2(t)=\rho^2(t) \text{ for } t\in [t_2,t_2+h],\\[2mm]
(x(t),y(t))=(x(t_2+h)-\mu (t-t_2-h)),y(t_2+h)) \text{ for }t\in [t_2+h, T].
\end{array}$$
 Then we have
$$\begin{array}{c}
x^2(T)+y^2(T)=\left( x(t_2)+\dot{x}(t_2)h+o(h)-\mu\left(\frac{\theta}{2}-h\right)\right)^2+(y(t_2)+\dot{y}(t_2)h+o(h))^2\\[2mm]
=
r^2_T+\Delta h+o(h)>r_T^2,\;\;\; h\in ]0, h_0].
\end{array}$$
This implies that $(x(T),y(T))\not\in C_T$. 

\vspace{5mm}

According to the necessary conditions given by Theorem \ref{main} there exist $\lambda\geq 0$, and $(p(\cdot),q(\cdot))\in BV$ such that $\lambda+|p(T)|+|q(T)|>0$, and
\begin{eqnarray}
\label{ex3}
&& dp=2\xi pdt+2\hat{x}d\eta,\\
&& dq=2\xi qdt+2\hat{y}d\eta,\label{ex4}\\
&& \max_{u\in[-\mu,1]} up=\hat{u}p,\label{ex5}\\
&& p(T)=\lambda-\beta \hat{x}(T),\;\;\; q(T)=-\beta\hat{y}(T),\;\;\;\beta\geq 0.\label{ex6}
\end{eqnarray}
Here $\xi\in L^{\infty}$, $\xi\geq 0$ and $\eta\in BV$. Observe that $\xi(t)=0$, $\eta(t)=({\rm const})$, if $\hat{r}(t)<\rho(t)$. 

Obviously $p(T)\leq 0$. Otherwise we have $\hat{u}(t)=1$, $t\in [T-\tilde{t},T]$ and with such  control the system cannot reach the terminal set from the initial condition. Notice that we consider $\mu$  small which  prevents  the system from reaching the boundary of $C(\cdot)$ with control $u(t)=-\mu$ in the initial instants of time.

 It is easy to see that any admissible trajectory satisfies the necessary conditions with $\lambda=1$, $p(\cdot)=0$, $q(t)=-\frac{y(T)}{x(T)}$,
$\xi(\cdot)=0$, and $\eta(\cdot)=0$. The only useful information that we retrieve from the necessary conditions is the equality $x^2(T)+y^2(T)=r^2_T$. Indeed, if $x^2(T)+y^2(T)<r^2_T$, then $\beta=0$, and (\ref{ex6}) implies $q(T)=0$ and $0=p(T)=\lambda$, a contradiction.
\end{example}

\begin{example}
Now, let us consider the same optimal control problem but with the dynamics governed by the following equations:
$$
\frac{d}{dt}\left(
\begin{array}{c}
x(t)\\
y(t)
\end{array}
\right)=
\left(
\begin{array}{c}
u(t)\\
-\sigma x(t)
\end{array}
\right)-N_{C(t)}(x(t),y(t)),
$$
where $\sigma>0$ is a small parameter.

 In this case the necessary conditions are: 
there exist $\lambda\geq 0$, $(p(\cdot),q(\cdot))\in BV$ such that $\lambda+|p(T)|+|q(T)|>0$, and
\begin{eqnarray}
\label{ex3a}
&& dp=\sigma qdt+2\xi pdt+2\hat{x}d\eta,\\
&& dq=2\xi qdt+2\hat{y}d\eta,\label{ex4a}\\
&& \max_{u\in[-\mu,1]} up=\hat{u}p,\label{ex5a}\\
&& p(T)=\lambda-\beta \hat{x}(T),\;\;\; q(T)=-\beta\hat{y}(T),\;\;\;\beta\geq 0.\label{ex6a}
\end{eqnarray}
Here $\xi\in L^{\infty}$, $\xi\geq 0$, and $\eta\in BV$. $\xi(t)=0$, with  $\eta(t)=({\rm const})$, if $\hat{r}(t)<\rho(t)$. Now the equality $p(t)=0$, 
$t\in [T-\tilde{t},T]$ is impossible. Indeed, from (\ref{ex3a}) we get $q(t)=0$, $t\in [T-\tilde{t},T]$. Then (\ref{ex6a}) implies that $\lambda=0$. Thus $p(t)<0$ and $\hat{u}(t)=-\mu$, $t\in [T-\tilde{t},T]$. In the interior of $C(\cdot)$, $p(\cdot)$ can change the sign only once. This means that $\hat{x}^2(T-\tilde{t})+\hat{y}^2(T-\tilde{t})=\rho^2(T-\tilde{t})$. In order to reach the point $(x(T-\tilde{t}),y(T-\tilde{t}))$ and to reach the boundary of $C(\cdot)$ from the initial point, the system must move with the control $\hat{u}=1$. This corresponds to $p(t)>0$. So, to find the optimal solution we have to analyze admissible  trajectories with the controls
$$
u(t)=\left\{
\begin{array}{cl}
1,& t\in [0,\tilde{t}],\\
-\mu, & t\in]\tilde{t},T],
\end{array}
\right.
$$
and choose the optimal value of $\tilde{t}$. Obviously all other controls give greater values of the functional.

\end{example}
\section*{Acknowledgements}

The authors gratefully acknowledge the financial support provided by the FEDER/COMPETE/NORTE2020/POCI/FCT funds 
[Grant number UID/EEA/00147/2019-SYSTEC] and 
the Portuguese Foundation for Science and Technology, in the framework of the Strategic Funding, 
through CFIS and CIDMA [Grant number UIDB/04106/2020].

\section*{Disclosure statement}

There are not declarations of interest.

\section*{Funding}

This work was supported by the FEDER/COMPETE/NORTE2020/POCI/FCT funds 
under Grant number [UID/EEA/00147/2019-SYSTEC] 
and 
the Portuguese Foundation for Science and Technology, in the framework of the Strategic Funding, 
under Grant number [UIDB/04106/2020].

%

%
%
%
%

\end{document}